\documentclass[12pt]{amsart}
\usepackage{amssymb}
\usepackage{amsthm}
\usepackage{euscript}
\usepackage{xspace}
\usepackage{enumerate}
\def\rmark{\mbox{$\rm\bf\rule{0.06em}{1.45ex}\kern-0.05em R$}}
\def\nmark{\mbox{$\rm\bf\rule{0.06em}{1.45ex}\kern-0.05em N$}}
\def\hmark{\mbox{$\rm\bf\rule{0.06em}{1.45ex}\kern-0.05em H$}}
\def\cmark{\mbox{$\rm\bf\kern0.2em\rule{0.06em}{1.45ex}\kern-0.3em C$}}

\theoremstyle{plain}
\newtheorem{lem}[subsection]{Lemma}

\newtheorem{thm}[subsection]{Theorem}
\newtheorem{prop}[subsection]{Proposition}

\newtheorem{cor}[subsection]{Corollary}

\begin{document}
\title[On modules of linear transformations]{On modules of linear transformations}
\author{M. Rahimi-Alangi   and Bamdad R. Yahaghi}
\address{ Department of Mathematics, Payame Noor University, P.O. Box 19395-3697 Tehran, Iran \newline
Department of Mathematics, Faculty of Sciences, Golestan University, Gorgan 19395-5746, Iran}
\email{mrahimi40@yahoo.com, \newline bamdad5@hotmail.com, bamdad@bamdadyahaghi.com}

\dedicatory{With kind regards,\\ dedicated to Rajendra Bhatia on the occasion of his sixtieth birthday}


\bibliographystyle{plain}

\begin{abstract}
Let $D$ be a division ring, $\mathcal V$ and $ \mathcal W$  vector
spaces over $D$, and ${\mathcal L(\mathcal V,\mathcal  W)}$ the
${\mathcal L(\mathcal W)}$-${\mathcal L(\mathcal V)}$  bimodule of all linear transformations from $\mathcal V$ into  $\mathcal W$. We prove some basic results about
 certain submodules of  $\mathcal L(\mathcal V, \mathcal W)$. For instance, we show, among other results, that a right submodule (resp. left submodule) of ${\mathcal L(\mathcal V,\mathcal  W)}$ is finitely generated whenever its image (resp. coimage) is finite-dimensional.
\end{abstract}
\maketitle \vspace{.5cm}

\bigskip

\begin{section}
{\bf One-sided submodules of linear transformations}
\end{section}

Throughout this note, unless otherwise stated, $D$ denotes a
division ring, $\mathcal V$ and $\mathcal W$ right (resp. left)
vector spaces over $D$, and $\mathcal L (\mathcal V, \mathcal W)$
the set of all right (resp. left) linear transformations $A:
{\mathcal V} \longrightarrow \mathcal W$ such that $A(x + y)= Ax +
Ay $ and $A (x\lambda) = (Ax) \lambda$ (resp. $A (\lambda x) =
\lambda (Ax)$) for all $ x , y \in \mathcal V$ and $\lambda \in D$.
From this point on, we only consider right vector spaces and
linear transformations among them  because everything remains valid if
one considers left vectors spaces and linear transformations acting among them.
It is well-known that $\mathcal L (\mathcal V, \mathcal W)$ forms
an abelian group under the addition of linear transformations.
When $ \mathcal V= \mathcal W$, we use the symbol $\mathcal
L(\mathcal V)$ to denote $\mathcal L (\mathcal V, \mathcal W)$.  It
is easy to see that the set $\mathcal L(\mathcal V)$ forms a ring
under the addition and multiplication of linear transformations
which are, respectively, defined by $(A+ B)(x) := Ax + Bx$ and
$(AB)(x) := A(Bx)$. It is also easily verified that
 $\mathcal L (\mathcal V, \mathcal W)$ is
a right  $\mathcal L (\mathcal V)$-module
(resp.   left $\mathcal L ( \mathcal W)$-module) via
the multiplication of linear transformations. Throughout, by
saying $\mathcal I$ is a right (resp. left) submodule of
 $\mathcal L (\mathcal V, \mathcal W)$, we mean $\mathcal I$
 is a right   $\mathcal L (\mathcal V)$-submodule
  (resp.   left $\mathcal L ( \mathcal W)$-submodule) of
  $\mathcal L (\mathcal V, \mathcal W)$.  By the {\it image} and the {\it kernel} of the family $\mathcal
F \subseteq \mathcal L (\mathcal V, \mathcal W)$, denoted by ${\rm im}(\mathcal F)$ and $\ker(\mathcal F)$, respectively,
we mean $\langle\{ Ax: A \in {{\mathcal F}}, x \in {{\mathcal V}}\}\rangle$ and
$\bigcap_{A \in {\mathcal F}} \ker A$. The {\it coimage} and {\it
cokernel} of the family $\mathcal F$, denoted by ${\rm coim}(\mathcal F)$
and ${\rm coker}(\mathcal F)$, respectively, are defined as ${\mathcal
V}/\ker \mathcal F$ and ${\mathcal W}/{\rm im}(\mathcal F)$. As is usual, we use the
symbol ${\mathcal V}'$ for $\mathcal L(\mathcal V, D)$. The members of $\mathcal V'$
are called linear functionals on $\mathcal V$. Also, when $\mathcal V$ is a
right (resp. left) vector space,  ${\mathcal V}'$ is a left (resp.
right) vector space over $D$ endowed with the addition and scalar
multiplication defined by $(f + g)(x) := f(x) + g(x)$ and $(\lambda
f) (x) := \lambda f(x)$ (resp. $(f \lambda ) (x) := f(x) \lambda$)
for all $ x \in \mathcal V$ and $\lambda \in D$.  {\it The second dual}
of $\mathcal V$, denoted by $\mathcal V''$, is the dual of $\mathcal V'$. The
space $\mathcal V''$ has the same chirality as that of $\mathcal V$ over
$D$. It is easily seen that $\mathcal V$ naturally imbeds into $\mathcal
V''$ via the natural mapping $\ \widehat{}  : \mathcal V \to {\mathcal V}''$
($x \mapsto \widehat{ x}$) defined by $\widehat{ x} (f)= f(x)$ for
all  $f \in \mathcal V'$, and that the natural mapping is an isomorphism
of the vector spaces $\mathcal V$ and $\mathcal V''$ if and only if the
space $\mathcal V$ is finite-dimensional. For a collection $\mathcal C$ of vectors in a
vector space $\mathcal V$ over $D$, $\langle {\mathcal C} \rangle$ is used to
denote {\it the linear subspace spanned by} $\mathcal
C$. For a subset $S$ of ${\mathcal V}$, we define $S^{\perp} := \{ f \in
{\mathcal V}' : f(S)=0\}$. It is plain that $S^{\perp}$ is a subspace of
${\mathcal V}'$. For $T \in \mathcal L(\mathcal V, \mathcal W)$, $T' \in \mathcal L(\mathcal W', \mathcal V')$
denotes the {\it adjoint} of $T$ which is defined by
$(T'f)(v):=f(Tv)$ where $ f \in {\mathcal W}', v \in \mathcal V$. For a subset $\mathcal{S}$ of $\mathcal L(\mathcal V, \mathcal W)$, it is not difficult to see that the map $ \phi: (\ker  \mathcal{S})^\perp \longrightarrow \left(\frac{\mathcal{V}}{\ker  \mathcal{S}}\right)'$ defined by $ \phi f(x+ \ker \mathcal{S}) = f(x)$, where $ f \in (\ker  \mathcal{S})^\perp$ and $ x \in \mathcal V$, is an isomorphism of vector spaces. Therefore,
$$(\ker  \mathcal{S})^\perp \cong \left(\frac{\mathcal{V}}{\ker  \mathcal{S}}\right)'.$$
By a weak right (resp. left) submodule of
$ \mathcal L(\mathcal W', \mathcal V')$, we mean $\mathcal I' := \{ T' \in \mathcal L(\mathcal W', \mathcal V') : T \in \mathcal I \}$,
where $\mathcal I$ is a left (resp. right) submodule of $\mathcal L (\mathcal V, \mathcal W)$. By definition
$\mathcal L' (\mathcal V, \mathcal W):= \{ T' \in  \mathcal L(\mathcal W', \mathcal V'): T \in \mathcal L (\mathcal V, \mathcal W)\}$.
An important subset of $\mathcal
L(\mathcal V, \mathcal W)$ is the class of rank-one linear
transformations. It can be shown that every rank-one linear
transformation in $\mathcal L(\mathcal V, \mathcal W)$ is of  the form $x \otimes f$ for some  $x \in
\mathcal W$ and $ f \in \mathcal V'$, where  $(x \otimes f)(y):= x
f(y)$ or $(x \otimes f)(y):= f(y)x$ depending on whether the space
$\mathcal W$ is a right or a left vector space over $D$. It is readily checked that $ (x \otimes f)' = f  \otimes  \widehat{ x}$.
Also, every finite-rank linear transformation is a finite sum of rank-one linear transformations. We use the symbol $\mathcal F (\mathcal V, \mathcal W)$ to denote the set, in fact the bi-module, of all finite-rank linear transformations from $\mathcal V$ into $\mathcal W$. As is usual,  $|A|$ is used to denote the cardinal number of the set $A$. We will make use of some basic results on cardinal numbers such as the Cantor--Schr$\ddot{{\rm o}}$der-Bernstein Theorem \cite[Theorem 1.5.3]{HSW}, that the class of cardinal numbers is well-ordered, that every two cardinal numbers are comparable, etc. We refer the reader to \cite{HJ} and \cite{HSW} for a general reference on set theory and cardinal arithmetic.
The following is a standard observation. If $ T \in \mathcal{L}(\mathcal{V}, \mathcal{W})$, then there are subsets
$B_1$ and $B_2$ of $\mathcal V$ such that $B_1$ and $B_1 \cup B_2$ are bases for $\ker T$ and $\mathcal{V}$, respectively. Moreover,  $T(B_2)$ is a basis for ${\rm im}(T)$ so that ${\rm im}( T)= \langle T(B_2)\rangle$.
 In fact, $\mathcal{V}=\langle B_1 \rangle \oplus \langle B_2 \rangle$. We refer the reader to \cite{H}, \cite{J1}, \cite{J2}, and \cite{R} for general references on  rings, modules, and linear algebra over division rings.

\bigskip

  Propositions \ref{1.1} and \ref{1.2} are likely known to the experts. For the counterpart of  Proposition \ref{1.2}(ii)-(iii) in the setting of Banach spaces, see  \cite{B}; also see \cite{D} and \cite{E}.

\bigskip

\begin{prop}\label{1.1}
{\rm  (i)}  {\it Let $S \in \mathcal{L}(\mathcal{V, W})$ and $T \in \mathcal{L}(\mathcal{V, Z})$. Then, $\ker S=\ker T$ iff there exists an injective linear transformation $P \in\mathcal{L}(\mathcal{W, Z})$  or $Q \in \mathcal{L}(\mathcal{Z, W})$ such that $T=PS$ or $S= QT$ depending on whether $\dim {\rm coker}(S) \leq \dim {\rm coker}(T)$ or $\dim {\rm coker}(T) \leq \dim {\rm coker}(S)$, respectively.}

{\rm (ii)}  {\it Let $S \in \mathcal{L}(\mathcal{V, W})$ and  $T \in \mathcal{L}(\mathcal{V, Z})$. Then, $\ker S \subseteq \ker T$ iff there exists a $P \in \mathcal{L}(\mathcal{W},\mathcal{Z})$ such that $T=PS$.}

{\rm (iii)}  {\it Let $S_i \in \mathcal{L}(\mathcal{V, W})$ ($1 \leq i \leq n$) and  $T \in \mathcal{L}(\mathcal{V, Z})$. Then, $\bigcap_{i=1}^n \ker S_i \subseteq \ker T$ iff there exist $P_i \in \mathcal{L}(\mathcal{W},\mathcal{Z})$ such that $T=P_1S_1 + \cdots + P_nS_n$. }
\end{prop}

\bigskip

\noindent {\bf Proof.} (i) The ``if" part is trivial. We prove the ``only if" part. Let $B_1$ be a basis for  $\ker S=\ker T$. Extend $B_1$ to a basis $ B_1 \cup B_2$
for $\mathcal{V}$. It follows that $\ker S=\langle B_1\rangle =\ker T , { \rm im}(S)=\langle S(B_2)\rangle$, and
$   {\rm im}(T)=\langle T(B_2)\rangle $. Extend   $S(B_2)$ and $T(B_2)$ to bases  $S(B_2) \cup B_3$ and
$T(B_2) \cup B_4 $ for $\mathcal W$ and  $\mathcal Z$, respectively. There are two cases to consider, namely $|B_3| \leq |B_4|$ or  $|B_4| \leq |B_3|$. If $|B_3| \leq |B_4|$, there is an injection $f: B_3 \rightarrow B_4$. Define the linear transformation
$P \in \mathcal{L (W, Z)}$ on $S(B_2) \cup B_3$ as follows
$$P(Sx) = Tx, P(y)= f(y), $$
for all $ x \in B_2$ and $ y \in B_3$, and extend $P$ to $\mathcal W$ linearly. It is plain that $PS=T$ and that
$P$ is injective because it takes the basis  $S(B_2) \cup B_3$ to a subset of the basis $T(B_2) \cup B_4 $  of $\mathcal W$. If $|B_4| \leq |B_3|$, the assertion follows in a similar fashion.

(ii) Again, the ``if" part is trivial. We prove the ``only if" part. Let $B_1$ be a basis for  $\ker S \subseteq \ker T$. Extend $B_1$ to a basis $ B_1 \cup B_2$ for
$ \ker T$ and then to a basis $ B_1 \cup B_2 \cup B_3$ for $\mathcal{V}$. It follows that $\ker S=\langle B_1\rangle$, $ \ker T= \langle B_1 \cup B_2\rangle$, ${ \rm im}(S)=\langle S(B_2 \cup B_3)\rangle$, and
$   {\rm im}(T)=\langle T(B_3)\rangle $. Extend   $S(B_2 \cup B_3)$ and $T(B_3)$ to bases
$S(B_2 \cup B_3) \cup B_4$ and
$T(B_3) \cup B_5 $ for $\mathcal W$ and $\mathcal{Z}$, respectively. Let $f: B_4 \rightarrow B_5$ be any function. Define the linear transformation
$P \in \mathcal{L (W, \mathcal{Z})}$ on $S(B_2 \cup B_3) \cup B_4$ as follows
$$P(Sx) = Tx, P(y)= f(y), $$
for all $ x \in B_2 \cup B_3$ and $ y \in B_4$, and extend $P$ to $\mathcal W$ linearly. It is plain that $PS=T$, as desired.

(iii) Just as in (ii), it suffices to prove the ``only if" part of the assertion. To this end, define $ S \in \mathcal{L}(\mathcal{V}, \mathcal{W}^n)$ and $T_1 \in \mathcal{L}(\mathcal{V}, \mathcal{Z}^n)$ on $\mathcal{V}$ by $Sx  = ( S_1 x, \ldots, S_n x)$ and $ T_1x = (Tx, \ldots, Tx)$, respectively. It is plain that $\ker S= \bigcap_{i=1}^n \ker S_i$ and $\ker T_1= \ker T$, which yields $\ker S \subseteq \ker T_1$. It thus follows from (ii) that there exists $ P \in \mathcal{L}(\mathcal{W}^n, \mathcal{Z}^n)$ such that $T_1 = PS$. If $P= (P_{ij})_{1 \leq i, j \leq n}$ is the standard matrix representation of $P$, where $P_{ij} \in \mathcal{L}(\mathcal{W},\mathcal{Z})$ ($1 \leq i, j \leq n$), we see that
$T = P_{11} S_1 + \cdots + P_{1n} S_n$, completing the proof.
\hfill\qed

\bigskip

\noindent {\bf Remark.}  In part (i) of the proposition, if $\dim {\rm coker}(S) = \dim {\rm coker}(T)$, then $\ker S=\ker T$ iff there exists an invertible  linear transformation $P \in\mathcal{L}(\mathcal{W, Z})$ such that $T=PS$.

 \bigskip

 \begin{prop}\label{1.2}
 {\rm (i)} {\it Let $S \in \mathcal{L}(\mathcal{V, W})$ and $T \in \mathcal{L}(\mathcal{Z, W})$. Then,  ${\rm im}(S)={\rm im}(T)$ iff there exists a  surjective  linear transformation $P \in\mathcal{L}(\mathcal{Z, V})$  or $Q \in \mathcal{L}(\mathcal{V, Z})$ such that $T=SP$ or $S= TQ$ depending on whether $\dim \ker S \leq \dim \ker T$ or  $\dim \ker T \leq \dim \ker S$, respectively. }

{\rm (ii)} {\it Let  $S \in \mathcal{L}(\mathcal{V, W})$ and $T \in \mathcal{L}(\mathcal{Z, W}) $. Then,  ${\rm im}(S)\subseteq {\rm im}(T)$ iff  there exists a  $P \in \mathcal{L}(\mathcal{V}, \mathcal{Z})$ such that $S=TP$.}

{\rm (iii)} {\it Let  $S \in \mathcal{L}(\mathcal{V, W})$ and $T_i \in \mathcal{L}(\mathcal{Z, W}) $ ($1 \leq i \leq n$). Then,  ${\rm im}(S)\subseteq {\rm im}(\{T_i\}_{i=1}^n)$ iff  there exist  $P_i \in \mathcal{L}(\mathcal{V}, \mathcal{Z})$ such that $S=T_1P_1 + \cdots +T_nP_n$.}

\end{prop}

  \bigskip

  \noindent {\bf Proof.} (i) The ``if" part being trivial, it suffices to prove the ``only if" part of the assertion. Let $B_1$ and $B_2$ be bases for  $\ker S$ and $\ker T$, respectively. Extend $B_1$ to a basis $ B_1 \cup B_3$ for $\mathcal{V}$. It follows that $\ker S=\langle B_1\rangle$,
$\ker T =  \langle B_2\rangle$,  ${ \rm im}(S)=\langle S(B_3)\rangle$. It is plain that $|S(B_3)| = |B_3|$. As
${\rm im}(S)={\rm im}(T)$, for each $ y \in B_3$, there exists an $ x_y \in \mathcal{V}$ such that
$Sy = Tx_y$. We see that $ B_4 = \{x_y \}_{y \in B_3}$ is linearly independent. To see this, suppose
 $x_{y_1} \lambda_1 + \cdots + x_{y_n}\lambda_n= 0$ for some $ n \in \mathbb{N}$, $ y_i \in B_3$,
  and $ \lambda_i \in D$, where $1 \leq i \leq n$. It follows that
  $Tx_{y_1} \lambda_1 + \cdots + Tx_{y_n}\lambda_n= 0$, from which we obtain
  $S(y_1 \lambda_1 + \cdots + y_n\lambda_n)= 0$, implying that
  $y_1 \lambda_1 + \cdots + y_n\lambda_n \in \ker S$. This yields $ \lambda_i = 0 $ for all $ 1 \leq i \leq n$
   because $ B_1 \cup B_3$ is a basis for $\mathcal{V}$ and $B_1$ is a basis for  $\ker S$.  We can
   also easily check that $B_2 \cup B_4$ is a basis for $\mathcal{Z}$. Now, there are two cases to consider, namely   $|B_1| \leq  |B_2|$ or   $|B_2| \leq  |B_1|$.  If $|B_1| \leq  |B_2|$,  there is a surjection  $f: B_2 \rightarrow B_1$. Define the linear transformation
   $P \in \mathcal{L (Z, V)}$ on $B_2 \cup B_4$   as follows
$$Px = f(x), Px_y= y, $$
for all $ x \in B_2$ and $ x_y \in B_4$, and extend $P$ to $\mathcal Z$ linearly. It is plain that $P$ is surjective
because it takes the basis  $B_2 \cup B_4$ of $\mathcal{Z}$ to the basis $B_1 \cup B_3 $  of $\mathcal V$. Also, we have $S=TP$ for
the following reason. That $ B_1 \cup B_3$ is a basis for
 $\mathcal{V}$.  That for all $ x \in B_1$, we have $ Px = f(x) \in B_2$, and hence $TPx= 0 = Sx$. And that for all  $y \in B_3$, we have  $Py= x_y$, which yields $TPy = Tx_y = Sy$. If $|B_2| \leq |B_1|$, the assertion follows in a similar fashion.

(ii) Again, it suffices to prove the ``only if" part of the assertion. Let $B_1$ and $B_2$ be bases for  $\ker S$ and $\ker T$, respectively. Extend $B_1$ to a basis $ B_1 \cup B_3$ for $\mathcal{V}$. It follows that $\ker S=\langle B_1\rangle$,
$\ker T =  \langle B_2\rangle$,  ${ \rm im}(S)=\langle S(B_3)\rangle$. It is plain that $|S(B_3)| = |B_3|$. As
${\rm im}(S) \subseteq {\rm im}(T)$, for each $ y \in B_3$, there exists an $ x_y \in \mathcal{Z}$ such that
$Sy = Tx_y$. Just as in (i), we see that $ B_4 = \{x_y \}_{y \in B_3}$ is linearly independent.
Let $f: B_1 \rightarrow B_2$ be any function. Define the linear transformation
   $P \in \mathcal{L (V, \mathcal{Z})}$ on $B_1 \cup B_3$   as follows
$$Px = f(x), Py= x_y, $$
for all $ x \in B_1$ and $ y \in B_3$, and extend $P$ to $\mathcal V$ linearly. We have $S=TP$ for
the following reason. That $ B_1 \cup B_3$ is a basis for
 $\mathcal{V}$.  That for all $ x \in B_1$, we have $ Px = f(x) \in B_2$, and hence $TPx= 0 = Sx$. And that for all  $y \in B_3$, we have  $Py= x_y$, which yields $TPy = Tx_y = Sy$.

 (iii) Just as in (ii), it suffices to prove the ``only if" part of the assertion. To this end, define $ S_1 \in \mathcal{L}(\mathcal{V}^n, \mathcal{W})$ and $T \in \mathcal{L}(\mathcal{Z}^n, \mathcal{W})$ by $S_1(x_1, \ldots, x_n)  =S x_1+  \cdots + S x_n$ and $ T(y_1, \ldots, y_n)  = T_1y_1+  \cdots + T_ny_n$, respectively. It is plain that ${\rm im}(S_1)= {\rm im}(S)$ and ${\rm im}(T)= {\rm im}(\{T_i\}_{i=1}^n)$, which yields ${\rm im}(S_1) \subseteq {\rm im}(T)$. It thus follows from (ii) that there exists $ P \in \mathcal{L}(\mathcal{V}^n, \mathcal{Z}^n)$ such that $S_1 = TP$. If $P= (P_{ij})_{1 \leq i, j \leq n}$ is the standard matrix representation of $P$, where $P_{ij} \in \mathcal{L}(\mathcal{V},\mathcal{Z})$ ($1 \leq i, j \leq n$), we easily see that
$S = T_1 P_{11} + \cdots + T_n P_{n1}$, completing the proof.
 \hfill\qed

\bigskip

\noindent {\bf Remark.} In part (i) of the proposition, if $\dim \ker S = \dim \ker T$, then ${\rm im}(S)={\rm im}(T)$ iff there exists an invertible  linear transformation $P \in\mathcal{L}(\mathcal{Z, V})$ such that $T=SP$.

 \bigskip

The proposition below is taken from \cite{RY2}, see \cite[Lemma 1.1]{RY2}. The counterpart of it in the setting of locally convex spaces was presented, with a detailed proof, in \cite{RY1}, see \cite[Lemma 1.1]{RY1}. It is however worth mentioning that \cite[Lemma 1.1]{RY2} was left as an exercise for the interested reader. The authors of  \cite{RY2} thought that the proof of the lemma is an imitation of that of \cite[Lemma 1.1]{RY1}. We have not been able yet to prove the second part of the proposition below for an arbitrary collection of linear transformations.

 \bigskip

\begin{prop}\label{1.3}
{\it Let $\mathcal V$ and $\mathcal W$ be  two vector spaces over
$D$ and ${\mathcal C}\subseteq \mathcal L(\mathcal V, \mathcal W)$. Then the following holds.

 {\rm (i)}  $ ({\rm im}\ {\mathcal C})^{\perp}= \ker({\mathcal C}'),$
where ${\mathcal C}'= \{ T': T \in {\mathcal C}\}$.
In other words,
$\langle \bigcup_{T \in {\mathcal C}} T {\mathcal V} \rangle^{\perp}= \bigcap_{T \in {\mathcal C}} \ker T'$.

 {\rm (ii)} If ${\mathcal C}= \{T_i\}_{i=1}^n$,  where $n \in \mathbb N$, then $(\ker {\mathcal C})^{\perp} = {\rm im}({\mathcal C}')$.
In other words, $(\bigcap_{i=1}^n \ker T_i)^{\perp} =
\langle \bigcup_{i=1}^n T_i' {\mathcal W}' \rangle$.
}
\end{prop}

\bigskip

\noindent {\bf Proof.} (i)  First, it is plain that $\langle \bigcup_{T \in {\mathcal C}} T {\mathcal V} \rangle^{\perp} \subseteq  \ker T'$ for each $T \in {\mathcal C}$, implying that $\langle \bigcup_{T \in {\mathcal C}} T {\mathcal V} \rangle^{\perp} \subseteq \bigcap_{T \in {\mathcal C}} \ker T'$. Now let $ g \in \bigcap_{T \in {\mathcal C}} \ker T'$ be arbitrary. It follows that $ T'g = 0$ for all $ T \in {\mathcal C}$. This yields  $ T'g({\mathcal V}) = g(T{\mathcal V})= 0$, implying that $ g \langle \bigcup_{T \in {\mathcal C}} T {\mathcal V} \rangle= 0$. That is, $ g \in \langle \bigcup_{T \in {\mathcal C}} T {\mathcal V} \rangle^{\perp}$. This proves the assertion.

(ii) First, it is plain that
$ (\bigcap_{T \in {\mathcal C}} \ker T)^{\perp} \supseteq T' {\mathcal W}'$,
for all $ T \in {\mathcal C}$, implying that
$ (\bigcap_{T \in {\mathcal C}} \ker T)^{\perp} \supseteq
\langle \bigcup_{T \in {\mathcal C}} T' {\mathcal W}' \rangle$. Now let $ f \in (\bigcap_{i=1}^n \ker T_i)^{\perp}$ be arbitrary. It follows that $\bigcap_{i=1}^n \ker T_i \subseteq \ker f$. It thus follows from Proposition \ref{1.1}(iii) that there are $ g_i \in {\mathcal W}'$ ($1 \leq i \leq n$) such that $ f = g_1 T_1 + \cdots + g_n T_n= T_1' g_1 + \cdots + T_n'g_n$, implying that $f \in \langle \bigcup_{i=1}^n T_i' {\mathcal W}' \rangle$. This completes the proof.
\hfill\qed

\bigskip

\noindent {\bf Remarks.} 1. It is easily shown that if $\mathcal V$ is finite-dimensional, then, using the reflexivity of $\mathcal V$, (ii) holds for all collections ${\mathcal C}\subseteq \mathcal L(\mathcal V, \mathcal W)$.

2. If $ n= \dim {\rm coim}({\mathcal C}) < \infty$, then there are $T_i \in \mathcal C$ ($1 \leq i \leq m \leq n$) such that $(\bigcap_{i=1}^m \ker T_i)^{\perp}= {\rm im}({\mathcal C}')$.

\bigskip

\begin{cor}\label{1.4}
{\rm (i)}  {\it Let $S_i \in \mathcal{L}(\mathcal{V, W})$ ($1 \leq i \leq n$) and  $T \in \mathcal{L}(\mathcal{V, Z})$. Then, $\bigcap_{i=1}^n \ker S_i \subseteq \ker T$  iff ${\rm im}(T')\subseteq {\rm im}(\{S_i'\}_{i=1}^n)$ iff there exist $P_i \in \mathcal{L}(\mathcal{W},\mathcal{Z})$ such that $T=P_1S_1 + \cdots + P_nS_n$.}

{\rm (ii)} {\it Let  $S \in \mathcal{L}(\mathcal{V, W})$ and $T_i \in \mathcal{L}(\mathcal{Z, W}) $ ($1 \leq i \leq n$). Then, ${\rm im}(S)\subseteq {\rm im}(\{T_i\}_{i=1}^n)$ iff $\bigcap_{i=1}^n  \ker T_i' \subseteq \ker S'$ iff there exist  $P_i \in \mathcal{L}(\mathcal{V}, \mathcal{Z})$ such that and $S=T_1P_1 + \cdots + T_nP_n$.}
\end{cor}

 \bigskip

\noindent {\bf Proof.} This is a quick consequence of Propositions \ref{1.1}(iii), \ref{1.2}(iii), and \ref{1.3}.
\hfill\qed

\bigskip

    \begin{thm}\label{1.5}
 {\it Let $\mathcal{V}$ and $\mathcal W$ be vector spaces over a division ring $D$ and $\mathcal I$ be a right submodule (resp.  left submodule) of  $\mathcal{L}(\mathcal{V, W})$. If  $\mathcal I$  finitely generated or  $\mathcal W$ (resp. $\mathcal{V}$) is finite-dimensional, then $\mathcal{I}= \{ T \in \mathcal{L}(\mathcal{V, W}): T \mathcal{V} \subseteq {\rm im}(\mathcal{I})\}$ (resp. $\mathcal{I}= \{ T \in \mathcal{L}(\mathcal{V, W}):  T\ker(\mathcal{I}) = \{0\} \}$). Moreover, if $\dim \mathcal{V} \geq \dim \mathcal W$ (resp. $\dim \mathcal{V} \leq \dim \mathcal W$), then every such right (res. left) submodule is principal. In particular, if $ \mathcal{V} = \mathcal W$, then every finitely generated one-sided ideal of $\mathcal{L}(\mathcal{V})$ is principal.
}
\end{thm}

 \bigskip

\noindent {\bf Proof.} First, let $\mathcal I$ be a right submodule of $\mathcal{L}(\mathcal{V, W})$. It follows from the hypothesis that there are $ T_i \in \mathcal I$ ($1 \leq i \leq m$) such that $ \mathcal{M} = {\rm im}(\mathcal{I})= {\rm im}(\{T_i\}_{i=1}^m)$. To prove the assertion, we need to show that if $ T \in \mathcal{L}(\mathcal{V, W})$ and $ T \mathcal{V} \subseteq \mathcal{M}$, then $T \in \mathcal{I}$. But this is a straightforward consequence of Proposition \ref{1.2}(iii), proving the assertion. If  $\mathcal W$ is finite-dimensional, then so is $ \mathcal{M} = {\rm im}(\mathcal{I})$ and hence $ \mathcal{M} = {\rm im}(\mathcal{I})= {\rm im}(\{T_i\}_{i=1}^m)$ for some $ T_i \in \mathcal I$ ($1 \leq i \leq m$). So the assertion follows from the above.

Next, let $\mathcal I$ be a left submodule of $\mathcal{L}(\mathcal{V, W})$. If $\mathcal I$ is finitely generated, then there are $ T_i \in \mathcal I$ ($1 \leq i \leq m$) such that $ \mathcal{M} = \ker(\mathcal{I})= \ker(\{T_i\}_{i=1}^m)$, in which case the assertion is a quick consequence of Proposition \ref{1.1}(iii). If $\mathcal{V}$ is finite-dimensional, then, in view of the remark following Proposition \ref{1.3}, we see that $(\ker {\mathcal I})^{\perp} = {\rm im}({\mathcal I}')$. We need to show that if $ T \in \mathcal{L}(\mathcal{V, W})$ and $ \mathcal{M}= \ker {\mathcal I} \subseteq \ker T$, then $T \in \mathcal{I}$. From $\ker {\mathcal I} \subseteq \ker T$, as $\mathcal{V}$ is finite-dimensional, taking perp of both sides of the inclusion, we obtain ${\rm im}(T') \subseteq {\rm im}({\mathcal I}')$. But ${\rm im}({\mathcal I}')$ is a subspace of the finite-dimensional space $\mathcal{V}'$. Consequently, there are $S_i \in \mathcal{I}$, such that  ${\rm im}({\mathcal I}') = {\rm im}(\{S_i'\}_{i=1}^n)$. Therefore, ${\rm im}(T') \subseteq {\rm im}(\{S_i'\}_{i=1}^n)$, and hence from Corollary \ref{1.4}, we obtain $P_i \in \mathcal{L}(\mathcal{W})$ ($1 \leq i \leq n$) such that $T=P_1S_1 + \cdots + P_nS_n$. This yields $T \in {\mathcal I}$, as desired.

The last part of the assertion readily follows from the first part of the assertion. Because  if $\dim \mathcal{V} \geq \dim \mathcal W$ (resp. $\dim \mathcal{V} \leq \dim \mathcal W$) and $\mathcal{M} = {\rm im}(\mathcal{I})$ (resp. $\mathcal{M} =\ker(\mathcal{I})$), there is always a $T_0 \in \mathcal{L}(\mathcal{V, W})$ such that ${\rm im} (T_0) = \mathcal{M}$ (resp. $\ker T_0 = \mathcal{M}$).
\hfill\qed

\bigskip

    \begin{thm}\label{1.6}
 {\it Let $\mathcal{V}$ and $\mathcal W$ be vector spaces over a division ring $D$ and $\mathcal I$ be a right submodule (resp.  left submodule) of  $\mathcal{L}(\mathcal{V, W})$ such that ${\rm im}(\mathcal{I})= {\rm im}(\{T_i\}_{i=1}^n)$ (resp. ${\rm
 ker}(\mathcal{I})= {\rm ker}(\{T_i\}_{i=1}^n) $), where $ n \in \mathbb N$ and $T_i \in \mathcal{I}$ ($1 \leq i \leq n$). Then $\mathcal{I}$ is generated by $ \{T_i\}_{i=1}^n$, and hence finitely generated. Therefore,  $\mathcal{I}= \{ T \in \mathcal{L}(\mathcal{V, W}): T \mathcal{V} \subseteq {\rm im}(\mathcal{I})\}$ (resp. $\mathcal{I}= \{ T \in \mathcal{L}(\mathcal{V, W}):  T\ker(\mathcal{I}) = \{0\} \}$). Moreover, if $\dim \mathcal{V} \geq \dim \mathcal W$ (resp. $\dim \mathcal{V} \leq \dim \mathcal W$), then $\mathcal I$ is principal.}
\end{thm}

 \bigskip

\noindent {\bf Proof.} First, let $\mathcal{I}$ be a right submodule  of
$\mathcal{L}(\mathcal{V, W})$ such that ${\rm im}(\mathcal{I})= {\rm im}(\{T_i\}_{i=1}^n)$ for some
$ n \in \mathbb N$ and $T_i \in \mathcal{I}$. Let $ S \in  \mathcal{I}$ be arbitrary. As ${\rm im}(S) \subseteq {\rm im}(\mathcal{I})= {\rm im}(\{T_i\}_{i=1}^n)$, from Proposition \ref{1.2}(iii), we see that there exist $P_i \in \mathcal{L}(\mathcal{V})$ such that $S=T_1P_1 + \cdots + T_n P_n$, proving the assertion.

Next, let $\mathcal  I$ be a left submodule  of
$\mathcal{L}(\mathcal{V, W})$ such that ${\rm ker}(\mathcal{I})= {\rm ker}(\{T_i\}_{i=1}^n)$ for some
$ n \in \mathbb N$ and $T_i \in \mathcal{I}$.  Let $ S \in  \mathcal{I}$ be arbitrary. As ${\rm ker}(\{T_i\}_{i=1}^n)= {\rm ker}(\mathcal{I})  \subseteq {\rm ker}(S) $, from Proposition \ref{1.1}(iii), we see that there exist $P_i \in \mathcal{L}(\mathcal{W})$ such that $S=P_1T_1 + \cdots + P_nT_n$, proving the assertion. The rest follows from Theorem \ref{1.5}.
\hfill\qed

 \bigskip

\noindent {\bf Remark.} If $\mathcal{V= W}$, the transformations $T_i \in \mathcal{I}$ in the theorem can be chosen to be idempotents. And such submodules, and in fact such one-sided ideals, of $\mathcal{L}(\mathcal{V})$, i.e., those one-sided ideals whose images or kernels are the same as those of a finite subset of the ideals, are principal.

 \bigskip

   \begin{lem}\label{1.7}
 {\it  Let $\mathcal{V}$ and $\mathcal W$ be vector spaces over a division ring $D$ and $\mathcal M \subseteq \mathcal V'$ a finite-dimensional subspace of $\mathcal V'$ such that $ \dim \mathcal W \geq \dim \mathcal M$.  Then, there exists a $ T \in \mathcal{L}(\mathcal{V, W})$ such that ${\rm im}(T') = \mathcal M$.}
\end{lem}

 \bigskip

\noindent {\bf Proof.} Let  $\mathcal{V}$ and $\mathcal W$ be right vector spaces and $\{f_i\}_{i=1}^n$ be a basis for $\mathcal M $. As $ \dim \mathcal W \geq \dim \mathcal M$, there is an independent subset $ \{y_i\}_{i=1}^n$ of  $\mathcal W$. Choose  $\{g_i\}_{i=1}^n$ in $\mathcal W'$ such that $g_i (y_j) = \delta_{ij}$, where  $ 1 \leq i , j \leq n$  and
$\delta$ denotes the Kronecker delta. Let $ T = y_1 \otimes f_1 + \cdots +  y_n \otimes f_n$. It is plain that $ T' =  f_1 \otimes \widehat{y_1} + \cdots +  f_n \otimes \widehat{y_n}$. We can write
$$  f_i \otimes \widehat{y_i} (g_j) =  \widehat{y_i} (g_j) f_i= g_j (y_i)f_i=  \delta_{ij}f_i,$$
for all $ 1 \leq i , j \leq n$. This yields $T'(g_i) = f_i$ for each  $ 1 \leq i  \leq n$, implying that $T'\mathcal W'= \langle f_i \rangle_{i=1}^n=\mathcal M $, as desired.
\hfill\qed

\bigskip

\begin{thm} \label{1.8}
 {\it  Let $\mathcal{V}$ and $\mathcal W$ be vector spaces over a division ring $D$ and $\mathcal I$ be a right submodule (resp.  left submodule) of  $\mathcal{L}(\mathcal{V, W})$.  If the image (resp. coimage) of $\mathcal{I}$ is finite-dimensional, then  $\mathcal I$ is finitely generated, and hence $\mathcal{I}= \{ T \in \mathcal{L}(\mathcal{V, W}): T \mathcal{V} \subseteq {\rm im}( \mathcal{I})\}$ (resp. $\mathcal{I}= \{ T \in \mathcal{L}(\mathcal{V, W}):  T\ker \mathcal{I} = \{0\} \}$). Moreover, if $\dim \mathcal{V} \geq \dim \mathcal W$ (resp. $\dim \mathcal{V} \leq \dim \mathcal W$), then $\mathcal{I}$ is principal.   }
\end{thm}

 \bigskip

 \noindent {\bf Proof.} First, suppose that $\mathcal{I}$ is a right  submodule of
 $\mathcal{L}(\mathcal{V, W})$. By the hypothesis, there are linearly independent vectors $y_1, \ldots, y_r$ in ${\rm im}( \mathcal{I})$
 such that ${\rm im}(\mathcal{I})=\langle y_1, \ldots,  y_r\rangle$, where $y_i=T_i(x_i)$ for some vectors $x_i\in \mathcal{V}$  and $ T_i \in \mathcal{I}$. It thus follows that $ {\rm im}( \mathcal{I})=  {\rm im}( \{T_i\}_{i=1}^r)$. So the assertion follows from Theorem \ref{1.6}.

 Next, suppose that $\mathcal{I}$ is a left  submodule of $\mathcal{L}(\mathcal{V, W})$. Set $\mathcal{I}' = \{ T' : T \in \mathcal{I}\}$. Then  $\mathcal{I}' $ is a weak right submodule of $\mathcal{L}(\mathcal{W', V'})$. By the proof of Proposition \ref{1.3}(ii), we have $ {\rm im}( \mathcal{I}') \subseteq (\ker  \mathcal{I})^\perp$. But as $ \dim {\rm coim}(  \mathcal{I})= \dim \frac{\mathcal{V}}{\ker  \mathcal{I}} $ is finite, we can write
 $$\dim (\ker  \mathcal{I})^\perp=   \dim \left(\frac{\mathcal{V}}{\ker  \mathcal{I}}\right)'=\dim \frac{\mathcal{V}}{\ker  \mathcal{I}}   < \infty.$$
 Consequently, $\dim  {\rm im}( \mathcal{I}') < \infty$, which implies $ {\rm im}( \mathcal{I}')=  {\rm im} \{ T'_i\}_{i=1}^n$ for some $ n \in \mathbb{N}$. We now prove that  $\mathcal{I}$ is finitely generated, proving the assertion. Now let $ T \in  \mathcal{I}$ be arbitrary. As $ {\rm im}(T') \subseteq {\rm im}( \mathcal{I}')=  {\rm im} \{ T'_i\}_{i=1}^n$, from Corollary \ref{1.4}(i), we see that $ T = P_1 T_1 + \cdots + P_n T_n$ for some $ P_i \in \mathcal{L}(\mathcal{W})$ ($1 \leq i \leq n$). This means  $  \mathcal{I} $ is generated by  $\{ T_i\}_{i=1}^n$, as desired. Finally, suppose $\dim \mathcal{V} \leq \dim \mathcal W$.  Thus $  \dim \mathcal W \geq \dim  {\rm im}( \mathcal{I}')$, and hence by  Lemma \ref{1.7}, there is a $T_0 \in \mathcal{L}(\mathcal{V, W})$ such that ${\rm im}(T_0') =  {\rm im} \{ T'_i\}_{i=1}^n$. It thus follows from Corollary \ref{1.4}(i) that there are $ P_i \in \mathcal{L}(\mathcal{W})$ ($1 \leq i \leq n$) such that $T_0 = P_1 T_1 + \cdots + P_n T_n$. But $\mathcal{I}$ is a left  submodule, so we obtain $T_0 \in \mathcal{I}$. Now let $ T \in  \mathcal{I}$ be arbitrary. As $ {\rm im}(T') \subseteq {\rm im}( \mathcal{I}')= {\rm im}(T_0')$, once again, from Corollary \ref{1.4}(i), we see that $ T = P_0  T_0$ for some $ P_0 \in \mathcal{L}(\mathcal{W})$. This means  $  \mathcal{I} =\mathcal{L}(\mathcal{W}) T_0 $. That is,  $\mathcal{I}$ is principal, finishing the proof.
 \hfill\qed

\bigskip

\begin{thm} \label{1.9}
 {\it  Let $\mathcal{V}$ and $\mathcal W$ be vector spaces over a division ring $D$ and $\mathcal I$ be a right submodule  of  $\mathcal{L}(\mathcal{V, W})$.  Then,  $\mathcal{I} \cap \mathcal{F}(\mathcal{V, W}) = \{ T \in \mathcal{F}(\mathcal{V, W}): {\rm im}(T) \subseteq {\rm im}( \mathcal{I})\}$.  }
\end{thm}

 \bigskip

\noindent {\bf Proof.} Let $\mathcal{I}$ be a right  submodule of
 $\mathcal{L}(\mathcal{V, W})$ and $ T \in \mathcal{F}(\mathcal{V, W})$ such that ${\rm im}(T) \subseteq {\rm im}( \mathcal{I})$. It suffices to show that $T \in \mathcal{I}$. As $T$ has finite rank, it follows that ${\rm im}(T) \subseteq {\rm im}( \{T_i\}_{i=1}^n)$
 for some $n \in \mathbb N$ and  $T_i \in \mathcal{I}$ ($1 \leq i \leq n$). It thus follows from Proposition \ref{1.2}(iii) that $ T = T_1 P_1 + \cdots + T_n P_n$ for some $ P_i \in  \mathcal{L}(\mathcal{V})$ ($1 \leq i \leq n$). Consequently, $T \in \mathcal{I}$, as desired.
 \hfill\qed

\bigskip

\noindent {\bf Remark.} We conjecture that the counterpart of the theorem holds for left submodules of $\mathcal{L}(\mathcal{V, W})$. That is, {\it if $\mathcal I$ is a left submodule of  $\mathcal{L}(\mathcal{V, W})$, then $\mathcal{I}\cap \mathcal{F}(\mathcal{V, W})= \{ T \in \mathcal{F}(\mathcal{V, W}):  \ker \mathcal{I} \subseteq \ker T \}$.}

\bigskip

In the proofs of the next two corollaries, we make use of some basic results on the arithmetic of cardinals.
The following is a slight generalization of the lemma presented on page 257 of \cite{J1} on $\mathcal{L}(\mathcal{V, W})$.
\bigskip

\begin{cor}\label{1.10}
 {\rm (i)} {\it Let $\mathcal{V, W, X, Y}$  be vector spaces over a division ring $D$ and $ S \in \mathcal{L}(\mathcal{V, W})$  and $ T \in \mathcal{L}(\mathcal{X, Y})$. Then, ${\rm rank} (S) \leq {\rm rank} (T)$ iff there exist $ P \in \mathcal{L}(\mathcal{Y, W})$ and $ Q \in \mathcal{L}(\mathcal{V, X})$ such that $S = PTQ$. Moreover, if $ S \in \mathcal{L}(\mathcal{V, W})$  and $ T_i \in \mathcal{L}(\mathcal{X, Y})$ ($1 \leq i \leq n$), then ${\rm rank} (S) \leq \sum_{i=1}^n {\rm rank} (T_i) $ iff there exist $ P_i \in \mathcal{L}(\mathcal{Y, W})$ and $ Q_i \in \mathcal{L}(\mathcal{V, X})$ such that $S = \sum_{i=1}^n P_iT_iQ_i$. }

{\rm (ii)} {\it Let $\mathcal{V}$ and $\mathcal W$ be vector spaces over a division ring $D$ and $ S, T \in \mathcal{L}(\mathcal{V, W})$. Then   ${\rm rank} (S) \leq {\rm rank} (T)$ iff there exist $ P \in \mathcal{L}(\mathcal{W})$ and  $ Q \in \mathcal{L}(\mathcal{V})$ such that $S = PTQ$. Moreover, if $ S, T_i \in \mathcal{L}(\mathcal{V, W})$ ($1 \leq i \leq n$), then, ${\rm rank} (S) \leq \sum_{i=1}^n {\rm rank} (T_i) $ iff there exist $ P_i \in \mathcal{L}(\mathcal{ W})$ and $ Q_i \in \mathcal{L}(\mathcal{V})$ such that $S = \sum_{i=1}^n P_iT_iQ_i$.}

\end{cor}

 \bigskip

\noindent {\bf Proof.} (i)  First, let ${\rm rank} (S) \leq {\rm rank} (T)$. Let $\{Sx_i\}_{i \in I}$ and $\{Ty_j\}_{j \in J}$ be bases for $ {\rm im}(S)$ and $ {\rm im}(T)$, respectively. It follows that $|I| \leq |J|$, and hence there is a surjection $f: J \rightarrow I$. Extend $\{Ty_j\}_{j \in J}$ to a basis $ \{Ty_j\}_{j \in J} \cup B$ for $\mathcal Y$. Define $ P \in \mathcal{L}(\mathcal{Y, W})$ as follows
$$P(Ty_j) = Sx_{f(j)}, \ P(y) = 0, $$
for all $ j \in J$ and $ y \in B$, and extend $P$ to $\mathcal{Y}$ linearly. As $f: J \rightarrow I$ is a surjection, we have $ {\rm im}(PT)={\rm im}(S)$. It thus follows from Proposition \ref{1.2}(ii) that there exists a linear transformation  $ Q \in \mathcal{L}(\mathcal{V, X})$ such that $S = PTQ$. Next, let $S = PTQ$ for some $ P \in \mathcal{L}(\mathcal{Y, W})$ and $ Q \in \mathcal{L}(\mathcal{V, X})$. Clearly, $ TQ\mathcal{V} \subseteq T  \mathcal{V}$, which obtains $ PTQ\mathcal{V} \subseteq PT  \mathcal{V}$. This implies  ${\rm rank} (S)= {\rm rank}(PTQ) \leq {\rm rank}(PT) \leq {\rm rank}(T)$, completing the proof.

Now, suppose ${\rm rank} (S) \leq \sum_{i=1}^n {\rm rank} (T_i) $. If there is a $1 \leq k \leq n$ such that ${\rm rank}(T_k) = \infty$, then by \cite[Corollary 1.5.12]{HSW} and \cite[Theorem 8.1.6]{HJ}, we obtain
$$\sum_{i=1}^n {\rm rank} (T_i)= \max_{1 \leq i \leq n} {\rm rank} (T_i)= {\rm rank} (T_j)$$
for some $ 1 \leq j \leq n$. Thus if ${\rm rank}(T_k) = \infty$ or ${\rm rank}(S) \leq {\rm rank}(T_k)$ for some $1 \leq k \leq n$, then the assertion follows from the proof above. So we may assume without loss of generality that $T_k$'s are all finite-rank linear transformations, and hence so is $S$, and that $ {\rm rank}(T_k) < {\rm rank}(S)$ for all $ 1 \leq k \leq n$. Let $\{Sx_i\}_{i \in I}$ and $\{T_ky_j\}_{j \in J_k}$ be bases for $ {\rm im}(S)$ and $ {\rm im}(T_k)$ ($1 \leq k \leq n$), respectively. It follows that $|I| \leq |J_1| + \cdots + |J_n|$ and $ |I| > |J_k|$ for each $k= 1, \ldots, n$. From this, as $I$ and $J_k$'s are all finite,  we easily see that there is an onto map $ f: \bigcup_{k=1}^n J_k \longrightarrow I$ such that $f|_{J_k}$ is one-to-one for each $k= 1, \ldots, n$. Extend $\{T_ky_j\}_{j \in J_k}$ to a basis $\{T_ky_j\}_{j \in J_k} \cup B_k$ for  $\mathcal Y$. Define $ P_k \in \mathcal{L}(\mathcal{Y, W})$ as follows
$$P_k(T_ky_j) = Sx_{f(j)}, \ P(y) = 0, $$
for all $ j \in J_k$ and $ y \in B_k$, and extend $P_k$ to $\mathcal{Y}$ linearly. It is plain that $ {\rm im} (S) = {\rm im} (P_1T_1, \ldots, P_nT_n)$. So the assertion follows from  Proposition \ref{1.2}(iii), as desired. The proof of the converse, which is quite straightforward, is omitted for brevity.

(ii) This is a quick consequence of (i).
\hfill\qed

\bigskip

\noindent {\bf Remark.} In view of the Rank-Nullity Theorem, it follows from part (ii) of the corollary that if $\mathcal{V}$ is finite-dimensional and $ S, T \in \mathcal{L}(\mathcal{V, W})$, then $\dim \ker S \geq \dim  \ker T$ iff there are  $ P \in \mathcal{L}(\mathcal{W})$ and $ Q \in \mathcal{L}(\mathcal{V})$ such that $S = PTQ$.

\bigskip

It is quite straightforward to check that every nonzero bi-submodule of $\mathcal{L}(\mathcal{V, W})$ includes all rank-one, and hence all finite-rank transformations. Therefore, if $\mathcal{V}$ or $\mathcal{W}$ is finite-dimensional, the trivial  bi-submodules of $\mathcal{L}(\mathcal{V, W})$, i.e.,  $0$ and $\mathcal{L}(\mathcal{V, W})$, are the only bi-submodules of $\mathcal{L}(\mathcal{V, W})$. In view of the preceding corollary, inspired by \cite[Theorem IX.5]{J1} or \cite[Theorem IV.17.1]{J2}, we state the following characterization of nontrivial  bi-submodules of $\mathcal{L}(\mathcal{V, W})$ whenever $\mathcal V$ and $\mathcal W$ are both infinite-dimensional.

\bigskip

\begin{cor}\label{1.11}
{\it Let $\mathcal{V}$ and $\mathcal W$ be infinite-dimensional vector spaces over a division ring $D$. Then the nontrivial  bi-submodules of $\mathcal{L}(\mathcal{V, W})$ are of the form
$$ \{ T \in \mathcal{L}(\mathcal{V, W}): {\rm rank} (T) < e\},$$
for some unique infinite cardinal number $ e \leq \min (\dim \mathcal{V}, \dim \mathcal{ W})$. }
\end{cor}

 \bigskip

\noindent {\bf Proof.}  Let $ \mathcal I$ be a nontrivial bi-submodule of $\mathcal{L}(\mathcal{V, W})$. Let $e_I$ be the smallest cardinal number such that $ e_I > {\rm rank}(A)$ for all $ A \in \mathcal I$. The existence of $e_I$ follows from the fact that the class of cardinal numbers is well-ordered, see \cite[Theorem 1]{Ho}, and its uniqueness follows from the Cantor-Schr$\ddot{{\rm o}}$der-Bernstein Theorem, see for instance \cite[Theorem 4.1.6]{HJ}. The cardinal number $e_I$ is infinite because $ \mathcal I$ includes all finite-rank linear transformations. If $ e_I > \min (\dim \mathcal{V}, \dim \mathcal{ W})$, then there is a $T \in \mathcal I$ such that $ {\rm rank}(T)= \min (\dim \mathcal{V}, \dim \mathcal{ W})$. This, in view of Corollary \ref{1.10}, would imply that $\mathcal I= \mathcal{L}(\mathcal{V, W})$, which is impossible. Thus, $ e_I \leq \min (\dim \mathcal{V}, \dim \mathcal{ W})$. Let $S \in \mathcal{L}(\mathcal{V, W})$ with ${\rm rank} (S) < e_I$ be arbitrary. We prove the assertion by showing that $S \in \mathcal I$. If ${\rm rank} (T) < \infty$, the assertion is trivial because  $ \mathcal I$ contains all finite-rank linear transformations. So we may assume that  $e_I >\aleph_0$, ${\rm rank} (S) < e_I$, and that $S$ has infinite rank. It follows from the definition of  $e_I$ that there is a $T \in \mathcal I$ such that ${\rm rank} (S) \leq {\rm rank} (T)$. So by Corollary \ref{1.10}, we have $ S= PTQ$ for some $ P \in \mathcal{L}(\mathcal{W})$ and  $ Q \in \mathcal{L}(\mathcal{V})$. Therefore, $S \in \mathcal I$, as desired. This completes the proof.
\hfill\qed

 \bigskip

\noindent  {\bf Acknowledgment.} The authors thank the referee for reading the manuscript carefully and for making helpful comments.

\end{document}